\documentclass[12pt]{amsart}

\usepackage{fullpage}
\usepackage{amsmath}
\usepackage{amssymb}
\usepackage{amsthm}
\usepackage{amstext}
\usepackage{appendix}
\usepackage{perpage}
\usepackage{color}
\usepackage{tikz-cd}
\usepackage[colorlinks=true,linkcolor=teal, citecolor=teal]{hyperref}

\newtheorem{thm}{Theorem}[section]
\newtheorem{lemm}[thm]{Lemma}
\newtheorem{coro}[thm]{Corollary}
\newtheorem{prop}[thm]{Proposition}

\theoremstyle{definition}

\newtheorem{defi}[thm]{Definition}

\begin{document}

\title{The algebraisation of higher Deligne--Lusztig representations}

\author{Zhe Chen \and Alexander Stasinski}

\address{Department of Mathematical Sciences, Durham University, Durham, DH1 3LE, UK}

\email{zhe.chen@durham.ac.uk}
\email{alexander.stasinski@durham.ac.uk}

\begin{abstract}
In this paper we study higher Deligne--Lusztig representations of reductive groups over finite quotients of discrete valuation rings. At even levels, we show that these geometrically constructed representations coincide with certain induced representations in the generic case; this gives a solution to a problem raised by Lusztig. In particular, we determine the dimensions of these representations. As an immediate application we verify a conjecture of Letellier for $\mathrm{GL}_2$ and $\mathrm{GL}_3$.
\end{abstract}

\maketitle

\tableofcontents

\section{Introduction}\label{sec:Intro}

In \cite{Lusztig1979SomeRemarks} Lusztig proposed a geometric (cohomological) construction (later proved in \cite{Lusztig2004RepsFinRings} for function fields and in \cite{Sta2009Unramified} in general) of representations of reductive groups over finite rings $\mathcal{O}_r=\mathcal{O}/\pi^r$, where $\mathcal{O}$ is the ring of integers in a non-archimedean local field with residue field $\mathbb{F}_q$, $\pi$ a uniformiser and $r\geq1$ a positive integer. This generalises the construction of Deligne and Lusztig \cite{DL1976} corresponding to the case $r=1$, which is the only known way to produce almost all irreducible representations of a general connected reductive group over a finite field. This generalised Deligne--Lusztig theory is a unified way to deal with all $r\geq1$. However, for $r>1$, besides the geometric construction, there is also a Clifford theoretic algebraic construction of representations of these groups. This algebraic method depends on the parity of $r$, and the idea can be traced back to Shintani \cite{Shintani1968sqr_int_irrep_lin_gr} and G\'erardin \cite{Gerardin1973GL_n}, who use this construction to study the representations of $p$-adic groups.

\vspace{2mm} Let $\mathbb{G}$ be a reductive group scheme over $\mathcal{O}_r$. For $r>1$, the geometrically constructed representations and the algebraically constructed representations share the same set of parameters, the pairs consisting of a maximal torus in $\mathbb{G}$, and a character of the $\mathcal{O}_r$-points of the torus satisfying certain regularity conditions (see Definition~\ref{def:regular} and \ref{def:gen-pos}). So a natural question, suggested by Lusztig in \cite[Section~1]{Lusztig2004RepsFinRings}, is whether the geometrically constructed representations coincide with the algebraically constructed representations. In Section~\ref{main result} we give a positive answer to this question for even levels $r=2l$. 

\vspace{2mm} Let $\mathcal{O}^{\mathrm{ur}}$ be the ring of integers in the maximal unramified extension of the field of fractions of $\mathcal{O}$, and put $\mathcal{O}^{\mathrm{ur}}_r=\mathcal{O}^{\mathrm{ur}}/\pi^r$. Denote the residue field of $\mathcal{O}^{\mathrm{ur}}$ by $k=\overline{\mathbb{F}}_q$. Let $\mathbf{G}$ be the base change of $\mathbb{G}$ to $\mathcal{O}^{\mathrm{ur}}_r$, and let $\mathcal{F}$ be the Greenberg functor from schemes of finite type over $\mathcal{O}^{\mathrm{ur}}_r$ to schemes over $k$. Then 
\begin{equation*}
G=G_r:=\mathcal{F}(\mathbf{G})
\end{equation*}
is a smooth affine algebraic group over $k$ such that $G(k)\cong\mathbf{G}(\mathcal{O}^{\mathrm{ur}}_r)$. Moreover, $G$ carries a Frobenius endomorphism $F$ such that
\begin{equation*}
G^F \cong \mathbb{G}(\mathcal{O}_r),
\end{equation*}
as finite groups. For a maximal torus in $\mathbb{G}$, we similarly obtain a subgroup $T$ of $G$. Throughout this paper we fix an arbitrary positive integer $r\geq 1$. For any integer $i$ such that $r\geq i\geq1$, let $\rho_{r,i}\colon G\rightarrow G_{i}$ be the reduction map modulo $\pi^i$; note that this is a surjective algebraic group morphism, and we denote the kernel by $G^{i}=G^i_r$. We also set $G^0=G$ (this is not the identity component $G^{\circ}$). Similar notation applies to closed subgroups of $G$.

\vspace{2mm} We now describe our main result. Let $\theta$ be a character of $T^F$. Assume that $r=2l$ is even. Then $G^l$ is abelian, and $T$ is a quotient of $TG^l$, so ${\theta}$ extends trivially to a character $\widetilde{\theta}$ of $(TG^l)^F$. Assume that $\theta$ is \emph{generic} (see Definition~\ref{generic}), our main result (see Theorem~\ref{multiplicity one, even levels} and Corollary~\ref{Main result at even levels}) says that
\begin{equation*}
R_{T}^{\theta} \cong \mathrm{Ind}_{(TG^{l})^F}^{G^F}\widetilde{\theta},
\end{equation*}
where $R_{T}^{\theta}$ is the higher Deligne--Lusztig representation (see Definition~\ref{Defi: DL rep}). As a consequence, $R_{T}^{\theta}$ has dimension $|G_l^F|/|T_l^F|$. The strategy of the proof is to first realise $\mathrm{Ind}_{(TG^{l})^F}^{G^F}\widetilde{\theta}$ as the cohomology of the Lang pre-image of certain algebraic group (see Proposition~\ref{algebraic realisation of arith rad rep}), and then show the inner product of these two representations equals $1$; the argument for the computation of inner product is generalised from the $\mathrm{GL}_n$ case in \cite{ZheChen_PhDthesis}. We remark that in the principal series case this isomorphism follows simply from the Mackey intertwining formula. We also remark that this isomorphism can fail when $\theta$ is not regular (this can be seen from the example computed by Lusztig in \cite[Section~3]{Lusztig2004RepsFinRings}). The case where $r$ is odd requires a different construction and is currently work in progress.

\vspace{2mm} Let $\mathfrak{g}$ be the Lie algebra of the reductive group $G_1$. For $\mathcal{O}_r=\mathbb{F}_q[[\pi]]/\pi^2$, by restricting the higher Deligne--Lusztig characters to the kernel $(G^1)^F\cong \mathfrak{g}^F$ one obtains invariant characters of finite Lie algebras. This was studied by Letellier in \cite{Let09CharRedLieAlg}, where he proposed several conjectures. One of them says roughly that any irreducible invariant character of $\mathfrak{g}^F$ appears in some Deligne--Lusztig character. We verify this conjecture for $\mathrm{GL}_2$ and $\mathrm{GL}_3$ in Section \ref{finite lie alg}. Previously this was only known for $\mathrm{GL}_2$ with the restriction that $|\mathbb{F}_q|>3$.

\vspace{2mm} During a summer school in Jul-Aug 2015, when we communicated with Lusztig about our methods and results, he told us that when he stated the expected relation between the algebraic and the geometric constructions, he had found a proof in the type $A_n$ case with $r=2$ (unpublished), by a method very different from ours.

\vspace{2mm} \noindent {\bf Acknowledgement.} We are grateful to A.-M.~Aubert and E.~Letellier for helpful discussions, and are thankful to G.~Lusztig for his encouragement. During the preparation of this work, {ZC} was partially supported by CSC/201308060137, and {AS} partially supported by EPSRC grant EP/K024779/1.

\section{Higher Deligne--Lusztig theory}\label{unrami DL}

Here we recall the main results developed in \cite{Lusztig2004RepsFinRings}, and \cite{Sta2009Unramified}. We preserve the notation introduced in Section~\ref{sec:Intro}: For $\mathbf{H}$ a smooth affine group scheme of finite type over $\mathcal{O}_r^{\mathrm{ur}}$, we have an associated algebraic group $H=H_r=\mathcal{F}\mathbf{H}$ over $k$, where $\mathcal{F}$ is the Greenberg functor; see \cite{Greenberg19611}, \cite{Greenberg19632}, \cite{Sta2009Unramified}, and \cite{Sta2012ReductiveGr} for its further properties. This $H$ is an affine smooth algebraic group over $k$ such that $H(k)\cong \mathbf{H}(\mathcal{O}^{\mathrm{ur}}_r)$.

\vspace{2mm} From now on, let $\mathbb{G}$ be a reductive group scheme over $\mathcal{O}_r$ (in other words, $\mathbb{G}$ is an affine smooth group scheme whose geometric fibre $\mathbf{G}_k$ is a connected reductive algebraic group in the classical sense; see e.g.\ \cite[XIX 2.7]{SGA3}). 

\vspace{2mm} Let $F\colon G\rightarrow G$ be a surjective algebraic group endomorphism such that the fixed points $G^F$ form a finite group; we call such a map a \emph{Frobenius endomorphism}. A closed subgroup $H\subseteq G$ is said to be $F$-rational (or rational when $F$ is fixed), if $F(H)\subseteq H$. In this paper we will only be concerned with the following typical situation: The Frobenius element $F$ in $\mathrm{Gal}(k/\mathbb{F}_q)$ extends to an automorphism of $\mathcal{O}^{\mathrm{ur}}_r$, and by the Greenberg functor this gives a rational structure on $G$ over $\mathbb{F}_q$. We denote the associated geometric Frobenius endomorphism again by $F$ and, as stated earlier, we thus have an isomorphism of finite groups $G^F\cong\mathbb{G}(\mathcal{O}_r)$. We write $L\colon g\mapsto g^{-1}F(g)$ for the Lang map associated to $F$.

\vspace{2mm} Let $\mathbf{T}\subset \mathbf{G}$ be a maximal torus such that $T=\mathcal{F}\mathbf{T}$ is $F$-rational, and let $\mathbf{B}$ be a Borel subgroup of $\mathbf{G}$ containing $\mathbf{T}$. Consider the Levi decomposition  $\mathbf{B}=\mathbf{U}\mathbf{T}$, where $\mathbf{U}$ is the unipotent radical of $\mathbf{B}$. The functor $\mathcal{F}$ gives a semi-direct product $B=\mathcal{F}\mathbf{B}=UT$ of closed subgroups of $G$, where $U=\mathcal{F}\mathbf{U}$. Let $\ell\neq p:=\mathrm{char}(\mathbb{F}_q)$ be a fixed prime number. We are interested in the \emph{higher level Deligne--Lusztig variety} associated to $T$ and $U$
\begin{equation*}
S_{T,U}=\{g\in G \mid g^{-1}F(g)\in FU\}=L^{-1}(FU),
\end{equation*}
where here, and in what follows, we often write $FU$ for $F(U)$. Note that $G^F\times T^F$ acts on $S_{T,U}$ by $(g,t)\colon x\mapsto gxt$, which induces an action on the compactly supported $\ell$-adic cohomology groups $H^i_c(S_{T,U}):=H^i_c(S_{T,U},\overline{\mathbb{Q}}_{\ell})$. 

\vspace{2mm} For any $\theta\in\widehat{T^F}=\mathrm{Hom}(T^F, \overline{\mathbb{Q}}_\ell^{\times})$, we denote by  $H^i_c(S_{T,U})_{\theta}$ the $\theta$-isotypical part of $H^i_c(S_{T,U})$. This is a $G^F$-submodule of $H^i_c(S_{T,U})$. We use the notation $H^*_c(-)$ for the alternating sum 
\begin{equation*}
H^*_c(-)=\sum_{i\geq 0}(-1)^iH^i_c(-,\overline{\mathbb{Q}}_{\ell}).
\end{equation*}
\begin{defi}\label{Defi: DL rep}
The \emph{higher Deligne--Lusztig representation} of $G^F$ associated to 
$\theta\in\widehat{T^F}$ is the virtual representation
\begin{equation*}
R^{\theta}_{T,U}=\sum_{i\geq 0}(-1)^iH_c^i(S_{T,U})_{\theta}.
\end{equation*}
In the situation we are interested, $R_{T,U}^{\theta}$ is independent of the choice of $U$; see Theorem~\ref{Irr DL}.
\end{defi}

The higher Deligne--Lusztig representations considered in this paper are the irreducible ones, or more precisely, the ones associated to certain characters of $T^F$ which are regular and in general position. We explain these notions.

\vspace{2mm} For any root $\alpha\in\Phi=\Phi(\mathbf{G},\mathbf{T})$ of 
$\mathbf{T}$, denote by $\mathbf{T}^{\alpha}$ the image of the coroot 
$\check{\alpha}$, and let $T^{\alpha}=\mathcal{F}\mathbf{T}^{\alpha}$. We 
write $\mathbf{U}_{\alpha}$ for the root subgroup of $\mathbf{U}$, and write 
$U_{\alpha}$ for its Greenberg functor image. For simplicity, we write $\mathcal{T}^{\alpha}$ for $(T^{\alpha})^{r-1}$. Note that $\mathbf{B}$ determines a subset $\Phi^{-}\subseteq\Phi$ of roots of $\mathbf{T}$. From now on we fix an arbitrary total order on $\Phi^-$.

\begin{defi}\label{def:regular}
Let $a$ be a fixed positive integer such that 
$F^a(\mathcal{T}^{\alpha})=\mathcal{T}^{\alpha}$ for every root $\alpha\in\Phi$ 
of $\mathbf{T}$. Consider the norm map $N^{F^a}_F(t):=t\cdot F(t)\cdots 
F^{a-1}(t)$ on $T^{F^a}$. Then $\theta\in\widehat{T^F}$ is called 
\emph{regular} if it is non-trivial on $N^{F^a}_F((\mathcal{T}^{\alpha})^{F^a})$ 
for every root $\alpha\in\Phi$. One knows that a regular character is regular 
with respect to any such $a$; see \cite[2.8]{Sta2009Unramified}.
\end{defi}

Since $\mathcal{O}^{\mathrm{ur}}_r$ is a strictly Henselian local ring, the reductive group scheme $\mathbf{G}$ is split with respect to every maximal torus (see \cite[2.1]{Sta2009Unramified}), therefore we can identify the Weyl group $W(T):=N(T)/T\cong W(T_1):=N(T_1)/T_1$; see \cite[XXII 3.4]{SGA3}.

\begin{defi}\label{def:gen-pos}
$\theta\in\widehat{T^F}$ is said to be in \emph{general position} if no non-trivial element in $W(T)^F=N(T)^F/T^F$ stabilises $\theta$.
\end{defi}

The following is one of the main results of \cite{Lusztig2004RepsFinRings} (in the function field case) and \cite{Sta2009Unramified} (in the general case).

\begin{thm}\label{Irr DL}
Suppose $\theta\in\widehat{T^F}$ is regular, then $R_{T,U}^{\theta}$ is independent of the choice of $U$, and if moreover $\theta$ is in general position, then $R^{\theta}_{T,U}$ is an irreducible representation up to sign.
\end{thm}
\begin{proof}
See \cite{Lusztig2004RepsFinRings} for the function fields and see \cite{Sta2009Unramified} for the general situation.
\end{proof}

\section{The algebraic construction}\label{arith radical}
From now on we assume $r=2l$ is even (note that $l$ is not the fixed prime 
$\ell$). Let $B_0=T_0U_0$ (resp.\ $T_0$, $U_0$) be the 
Greenberg functor image of a Borel subgroup $\mathbf{B}_0$ (resp.\ maximal 
torus 
$\mathbf{T}_0$, unipotent radical $\mathbf{U}_0$) of $\mathbf{B}$, such that 
$B_0$ is $F$-rational.
 Let ${\lambda}\in G$ be such that  $B={\lambda}B_0{\lambda}^{-1}$ and 
$T={\lambda}T_0{\lambda}^{-1}$. Note that 
${\lambda}^{-1}F({\lambda})=\hat{w}\in N(T_0)$ is a lift of some Weyl element 
$w\in W(T_0)$.

\begin{defi}
Along with the above notations, we denote by $U^{\pm}$ the commutative unipotent group $(U^{-})^lU^l$, and call it the \emph{arithmetic radical} 
associated to $T$. 
\end{defi}
Note that $T=\mathcal{F}\mathbf{T}$ is 
usually not a torus, but we sometimes still call it a torus. For 
convenience, we similarly say ``Borel subgroup'' for
$B=\mathcal{F}\mathbf{B}$.

\begin{lemm}\label{rationality of arith rad}
$U^{\pm}$ is normalised by $N(T)$, and it is $F$-rational.
\end{lemm}
\begin{proof}
Note that $U^{\pm}=\prod_{\alpha\in\Phi}{(U_{\alpha})^l}$, where 
$U_{\alpha}=\mathcal{F}\mathbf{U}_{\alpha}$ is defined above. For any $v\in 
W(T)$ we have $\hat{v}U_{\alpha}\hat{v}^{-1}=U_{v(\alpha)}$, where $\hat{v}$ is 
a lift of $v$ in $N(T)$, so 
\begin{equation*}
\hat{v}\prod_{\alpha\in\Phi}{(U_{\alpha})^l}\hat{v}^{-1}=\prod_{\alpha\in\Phi}{
(U_{v(\alpha)})^l}=\prod_{\alpha\in\Phi}{(U_{\alpha})^l},
\end{equation*}
which means $U^{\pm}$ is normalised by $N(T)$. Similarly, 
\begin{equation*}
F(U^{\pm})=F({\lambda}U^{\pm}_0{\lambda}^{-1})={\lambda}\hat{w}\prod_{
\beta\in\Phi_0}(U_{0,\beta})^l\hat{w}^{-1}{\lambda}^{-1},
\end{equation*}
where $\Phi_0$ is the root system for $\mathbf{T}_0$. The right hand side is 
\begin{equation*}
{\lambda}\prod_{\beta\in\Phi_0}(U_{0,w(\beta)})^l{\lambda}^{-1}={\lambda}\prod_{
\beta\in\Phi_0}(U_{0,\beta})^l{\lambda}^{-1}=U^{\pm}.
\end{equation*}
This proves the rationality.
\end{proof}

The variety $L^{-1}(U^{\pm})$ admits a left $G^F$-action and a right $T^F$-action, so $H^*_c(L^{-1}(U^{\pm}))$ is a $G^F\times T^F$-module. 

\begin{prop}\label{algebraic realisation of arith rad rep}
For every $\theta\in\widehat{T^F}$ we have 
$H^*_c(L^{-1}(U^{\pm}))_{\theta}\cong\mathrm{Ind}_{(TU^{\pm})^F}^{G^F}\widetilde
{\theta}$, where $\widetilde{\theta}$ is the trivial lift of $\theta$ from 
$T^F$ 
to $(TU^{\pm})^F$ (that is, $\widetilde{\theta}$ is the pull-back of 
$(TU^{\pm})^F\rightarrow T^F$).
\end{prop}
\begin{proof}
This is an argument analogous to the last paragraph in \cite[p.~81]{DM1991}. 
Consider the natural morphism $L^{-1}(U^{\pm})\rightarrow G/U^{\pm}$ given by 
$g\mapsto gU^{\pm}$. Note that $F(gU^{\pm})=gL(g)U^{\pm}=gU^{\pm}$, so its 
image 
is $(G/U^{\pm})^F\cong G^F/(U^{\pm})^F$. Note that its fibres are isomorphic to 
an affine space ($\cong U^{\pm}$), therefore 
$H^*_c(L^{-1}(U^{\pm}))\cong\overline{\mathbb{Q}}_{\ell}[G^F/(U^{\pm})^F]$ by 
basic properties of $\ell$-adic cohomology. Finally, 
$\overline{\mathbb{Q}}_{\ell}[G^F/(U^{\pm})^F]\otimes_{\overline{\mathbb{Q}}_{
\ell}[T^F]}\theta\cong\overline{\mathbb{Q}}_{\ell}[G^F]\otimes_{\overline{
\mathbb{Q}}_{\ell}[(TU^{\pm})^F]}\widetilde{\theta}$ as 
$\overline{\mathbb{Q}}_{\ell}[G^F]$-modules, thus 
$H^*_c(L^{-1}(U^{\pm}))_{\theta}\cong\mathrm{Ind}_{(TU^{\pm})^F}^{G^F}\widetilde
{\theta}$.
\end{proof} 

The representations $\mathrm{Ind}_{(TU^{\pm})^F}^{G^F}\widetilde{\theta}$ were already considered by G\'erardin in a more restrictive situation (i.e.\ $\mathbb{G}$ is defined over the field of fractions of $\mathcal{O}$, and it is split and its derived subgroup is simply-connected, and moreover, he required the maximal tori to be ``special'' in the sense of \cite[3.3.9]{Gerardin1975SeriesDiscretes}); see \cite{Gerardin1973GL_n} and \cite{Gerardin1975SeriesDiscretes}.

\vspace{2mm} Note that since $(G^l)^F$ is abelian, one has $(TU^{\pm})^F\subseteq\mathrm{Stab}_{G^F}(\widetilde{\theta}|_{(G^l)^F})$. By Clifford theory, if the equality holds, then $\mathrm{Ind}_{(TU^{\pm})^F}^{G^F}\widetilde{\theta}$ is irreducible; under G\'erardin's conditions the irreducibility always holds (see \cite[4.4.1 and 4.4.6]{Gerardin1975SeriesDiscretes}). We combine this property into the below definition.

\begin{defi}\label{generic}
A character $\theta\in\widehat{T^F}$ is \emph{generic} if it is regular, in general position, and $\mathrm{Stab}_{G^F}(\widetilde{\theta}|_{(G^l)^F})=(TU^{\pm})^F$.
\end{defi}

Actually, in some situation the stabiliser condition $(TU^{\pm})^F=\mathrm{Stab}_{G^F}(\widetilde{\theta}|_{(G^l)^F})$ is equivalent to the regularity of $\theta$, and implies the general position condition; we verify this for the Coxeter torus in a general linear group.

\begin{prop}\label{coxeter generic}
For $\mathbb{G}=\mathrm{GL}_n$ over $\mathcal{O}_r$, let $\mathbf{T}\subset\mathbf{G}$ be a maximal torus corresponding to the Coxeter element $w=(1,2,...,n)$. Then for $\theta\in\widehat{T^F}$,  the stabiliser condition $(TU^{\pm})^F=\mathrm{Stab}_{G^F}(\widetilde{\theta}|_{(G^l)^F})$ is equivalent to the regularity of $\theta$, and they imply $\theta$ is in general position.
\end{prop}
\begin{proof}
Note that $(G^l)^F\cong M_n(\mathcal{O}_l)$ (in the below we always assume 
this identification), and its irreducible characters are of the form 
$\psi_{\beta}(-)=\psi(\mathrm{Tr}(\beta(-)))$, where $\beta\in 
M_n(\mathcal{O}_l)$, and $\psi$ is a fixed complex-valued additive character 
on $\mathcal{O}_l$ which is non-trivial on the ideal $(\pi^{l-1})$. 
Suppose $\widetilde{\theta}|_{(G^l)^F}=\psi_{\beta}$, then the condition 
$(TU^{\pm})^F=\mathrm{Stab}_{G^F}(\widetilde{\theta}|_{(G^l)^F})$ is equivalent 
to $C_{\mathbb{G}(\mathcal{O}_l)}(\beta)=T_l^F$.

\vspace{2mm} Since $\widetilde{\theta}$ is trivial on $(U^{\pm})^F$, we have 
$\beta\in T_l^F$, so we can assume
\begin{equation*}
\beta_0=\lambda^{-1}\beta\lambda=\mathrm{diag}({\beta_1,...,\beta_n})\in 
M_n(\mathcal{O}^{\mathrm{ur}}_l)
\end{equation*}
(here the image of ${\lambda}$ modulo $\pi^l$ is again denoted by $\lambda$; 
this should make no confusion). With these notations the condition 
$C_{\mathbb{G}(\mathcal{O}_l)}(\beta)=T_l^F$ is equivalent to 
$C_{\lambda^{-1}\mathbb{G}(\mathcal{O}_l)\lambda}(\beta_0)=\lambda^{-1}
T_l^F\lambda$. However, as $\lambda^{-1}T_l^F\lambda\subseteq T_0$, this 
happens if and only if  $\beta_i-\beta_j$ is invertible for any $i\neq j$, and 
in particular $\theta$ is in general position.

\vspace{2mm} As we are concerning general linear groups, we can assume 
${\lambda}$ satisfies ${\lambda}^{-1}F({\lambda})=\hat{w}\in N(T)^F$; denote by 
$v$ the image of $\hat{w}$ in $\mathrm{GL}_n(\mathcal{O}^{\mathrm{ur}}_l)$. 
For any $t\in T_l$, we have $F(t)=\lambda 
vF(\lambda^{-1}t\lambda)v^{-1}\lambda^{-1}$. Denote by $F'$ the endomorphism 
$F'(g)=vF(g)v^{-1}$, then for any root $\alpha$, and any positive 
integer $m$ such that $F^m(\mathcal{T}^{\alpha})=\mathcal{T}^{\alpha}$, we have
\begin{equation*}
N_{F}^{F^m}(t)=\lambda t_0F'(t_0)\cdots F'^{m-1}(t_0)\lambda^{-1},
\end{equation*}
where $t\in (\mathcal{T}^{\alpha})^{F^m}$ and $t_0=\lambda^{-1}t\lambda$. Thus 
the regularity of $\theta$ is equivalent to:
\begin{equation*}
\psi(\mathrm{Tr}(\beta_0N^{F'^m}_{F'}(t_0)))\neq1
\end{equation*}
for some $t_0\in (\lambda^{-1}\mathcal{T}^{\alpha}\lambda)^{F'^m}$ for each 
given root $\alpha$ and $m$.

\vspace{2mm} Note that any conjugation from $T$ to $T_0$ takes the ``root 
subgroup'' $T^{\alpha}$ to a ``root subgroup'' of $T_0$, so we can write 
$t_0=\mathrm{diag}(0,...0,x,0,...0,-x,0,...,0)\in 
M_n(\mathcal{O}^{\mathrm{ur}}_l)$ for any 
$t_0\in\lambda^{-1}\mathcal{T}^{\alpha}\lambda$, where 
$x\in\pi^{l-1}\mathcal{O}^{\mathrm{ur}}_l\cong k$ is at position $(a,a)$ 
and $-x$ is at position $(b,b)$. As $v$ is a Coxeter element, we can take 
$m=n$, thus
\begin{equation*}
\mathrm{Tr}(\beta_0N^{F'^m}_{F'}(t_0))=\sum_{l=0}^{n-1}(\beta_{v^l(a)}-\beta_{v^
{l}(b)})F^l(x).
\end{equation*}
Since we are concerning the Coxeter element $(1,...,n)$, we can write 
$\beta_1=\beta'\in (\mathcal{O}^{\mathrm{ur}}_l)^{F^n}$ and $\beta_i=F^{i-1}(\beta')$; this enables us 
to rewrite the above as
\begin{equation*}
\mathrm{Tr}(\beta_0N^{F'^m}_{F'}(t_0))=\sum_{l=0}^{n-1}F^l(F^a(\beta'-F^{b-a}
(\beta'))x).
\end{equation*}
Therefore the regularity is equivalent to that, for any $b-a\in [1,...,n-1]$, the element
$\beta'-F^{b-a}(\beta')\in \mathcal{O}^{\mathrm{ur}}_l$ is invertible, i.e. 
$\beta_i-\beta_j$ is invertible for all $i\neq j$, and we see from the above 
this is equivalent to the stabiliser condition.
\end{proof}

\section{The main result}\label{main result}

\vspace{2mm} As before, $\mathbb{G}$ is a reductive group scheme over $\mathcal{O}_r$, $F$ is the corresponding Frobenius on $G$ and $\mathbf{T}$ is a maximal torus in $\mathbf{G}$ such that $T$ is $F$-rational. Moreover, $U$ is the Greenberg functor image of the unipotent radical of a Borel subgroup $\mathbf{B}$ of $\mathbf{G}$ containing $\mathbf{T}$. For any $v\in W(T)$, we fix a lift $\hat{v}\in N(T)$. Recall that (see Lemma \ref{rationality of arith rad}) $F(U^{\pm})=U^{\pm}$ and $\hat{v}U^{\pm}\hat{v}^{-1}=U^{\pm}$. Given two elements $x$ and $y$ in a group, we sometimes use the shorthand notation $x^y:=y^{-1}xy$ and ${^yx}:=yxy^{-1}$ for conjugations.

\vspace{2mm} Now we are going to present our main result. We start with the computation of inner products of Deligne--Lusztig representations and the representations produced from the arithmetic radicals.

\begin{thm}\label{multiplicity one, even levels}
Suppose $r=2l$ is even and $\theta\in\widehat{T^F}$ is regular and in general position. Then
\begin{equation*}
\langle \mathrm{Ind}_{(TU^{\pm})^F}^{G^F}\widetilde
{\theta},R_{T,U}^{\theta}\rangle_{G^F}=1.
\end{equation*}
\end{thm}
\begin{proof}
We want to compare the cohomology of $S_{T,U}=L^{-1}(FU)$ with the cohomology of the Lang pre-image $L^{-1}(FU^{\pm})$ of the arithmetic radical (see Proposition~\ref{algebraic realisation of arith rad rep}). One has
\begin{equation*}
\langle H^*_c(L^{-1}(FU^{\pm}))_{\theta},R_{T,U}^{\theta}\rangle_{G^F}=\dim 
H^*_c(\Sigma)_{\theta^{-1},\theta},
\end{equation*}
where 
\begin{equation*}
\Sigma:=\{(x,x',y)\in U^{\pm}\times FU\times G\mid xF(y)=yx'\}.
\end{equation*}
This follows from the $T^F\times T^F$-equivariant isomorphism
\begin{equation*}
G^F\backslash L^{-1}(U^{\pm})\times L^{-1}(FU)\cong \Sigma,\quad  
(g,g')\mapsto (g^{-1}F(g),{g'}^{-1}F(g'),g^{-1}g')
\end{equation*}
and the K\"{u}nneth formula; here $T^F\times T^F$ acts on $\Sigma$ by $(t,t')\colon (x,x',y)\mapsto (x^t,(x')^{t'},t^{-1}yt')$.

\vspace{2mm} The Bruhat decomposition $G_1=\coprod_{v\in W(T)} B_{1}\hat{v}B_{1}$ of $G_1=\mathbf{G}(k)$ gives the finite stratification (see, 
e.g.\ the proof of \cite[Lemma 2.3]{Sta2009Unramified}) $G=\coprod_{v\in 
W(T)}G_v$, where 
\begin{equation*}
G_v:=(U\cap\hat{v}U^{-}\hat{v}^{-1})(\hat{v}(U^{-})^1\hat{v}^{-1})\hat{v}TU,
\end{equation*}
and hence a finite partition into disjoint locally closed subvarieties
\begin{equation*}
\Sigma=\coprod_{v\in W(T)} \Sigma_{v},
\end{equation*}
where 
\begin{equation*}
\Sigma_{v}:=\{(x,x',y)\in U^{\pm}\times FU\times G_v\mid xF(y)=yx'\}.
\end{equation*}
For each $v$, consider the variety
\begin{equation*}
\mathcal{Z}_v:=(U\cap \hat{v}U^{-}\hat{v}^{-1})\times 
\hat{v}(U^{-})^{1}\hat{v}^{-1};
\end{equation*}
this allows us to consider
\begin{equation*}
\widehat{\Sigma}_v:=\{(x,x',u',u^{-},\tau,u)\in U^{\pm}\times FU\times \mathcal{Z}_v \times T\times U\mid xF(u'u^{-}\hat{v}\tau u)=u'u^{-}\hat{v}\tau 
ux'\}.
\end{equation*}
This is a locally trivial fibration $\widehat{\Sigma}_v\rightarrow\Sigma_v$ by an affine space ($\cong U\cap\hat{v}(U^-)^{1}\hat{v}^{-1}$), on which $T^F\times T^F$ acts as
\begin{equation*}
(t,t')\colon (x,x',u',u^{-},\tau,u)\longmapsto ({t}^{-1}xt,{t'}^{-1}x't',t^{-1}u't,t^{-1}u^{-}t,(t^{\hat{v}})^{-1}\tau t',{t'}^{-1}ut').
\end{equation*}
By the change of variable $x'F(u)^{-1}\mapsto x'$ we can rewrite $\widehat{\Sigma}_v$ as
\begin{equation*}
\widehat{\Sigma}_v=\{(x,x',u',u^{-},\tau,u)\in U^{\pm}\times FU\times \mathcal{Z}_v\times T\times U\mid xF(u'u^{-}\hat{v}\tau)=u'u^{-}\hat{v}\tau ux'\},
\end{equation*}
on which the $T^F\times T^F$-action does not change.

\vspace{2mm} For $i=0,1,...,r-1$ let $\mathcal{Z}_v(i)$ be the pre-image of $(\hat{v}U^{-}\hat{v}^{-1})^{i}=\hat{v}(U^{-})^{i}\hat{v}^{-1}$ under the product morphism 
\begin{equation*}
\mathcal{Z}_v=(U\cap \hat{v}U^{-}\hat{v}^{-1})\times 
\hat{v}(U^{-})^{1}\hat{v}^{-1}\longrightarrow \hat{v}U^{-}\hat{v}^{-1}.
\end{equation*}
Recall that for $i=0$ we always let $G^0=G$ for an algebraic group $G$. For each $v$ consider the partition $\widehat{\Sigma}_v=\Sigma'_v\sqcup\Sigma''_v$ of locally closed subvarieties, where 
\begin{equation*}
\Sigma'_v:=\{(x,x',u',u^{-},\tau,u)\in \widehat{\Sigma}_v\mid (u',u^{-})\in 
\mathcal{Z}_v\setminus\mathcal{Z}_v(l) \}
\end{equation*}
and
\begin{equation*}
\Sigma''_v:=\{(x,x',u',u^{-},\tau,u)\in \widehat{\Sigma}_v\mid (u',u^{-})\in 
\mathcal{Z}_v(l) \}.
\end{equation*}
Our goal is to show:
\begin{itemize}
\item[(a)] $\dim H^*_c(\Sigma''_v)_{\theta^{-1},\theta}=\left\{
\begin{array}{ll}
1 & v=1\\
0 & v\neq1,
\end{array}
\right.$

\item[(b)] $\dim H^*_c(\Sigma'_v)_{\theta^{-1},\theta}=0$ for all $v$.
\end{itemize}

\vspace{2mm} We start with (a), which is much easier.

\begin{lemm}\label{lemm a}
(a) is true.
\end{lemm}
\begin{proof}
Note that for any $(x,x',u',u^{-},\tau,u)\in\Sigma''_v$ we have
\begin{equation*}
u'u^{-}\in\hat{v}(U^{-})^{l}\hat{v}^{-1}\subseteq U^{\pm}=FU^{\pm},
\end{equation*}
so we can apply the changes of variables $(u'u^{-})^{-1}x\mapsto x$, and then 
$xF(u'u^{-})\mapsto x$. This allows us to rewrite $\Sigma''_v$ as 
\begin{equation*}
{\widetilde{\Sigma}''}_v :=\{(x,x',u',u^{-},\tau,u)\in U^{\pm}\times FU\times 
\mathcal{Z}_v(l)\times T\times U\mid xF(\hat{v}\tau)=\hat{v}\tau ux'\},
\end{equation*}
on which $T^F\times T^F$ acts in the same way as before. 

\vspace{2mm} Consider the algebraic group 
\begin{equation*}
H=\{(t,t')\in T_1\times T_1\mid tF(t^{-1})=F(\hat{v})t'F(t')^{-1}F(\hat{v}^{-1})\}.
\end{equation*}
Note that the action of $T_1^F\times T_1^F$ on $\widetilde{\Sigma}''_v$ extends to an action of $H$ (the torus $T_1$ is always a subgroup of $T$) in a natural way. The identity component $H^{\circ}$ is a torus acting on ${\widetilde{\Sigma}''}_v$, and thus by basic properties of $\ell$-adic cohomology we have 
\begin{equation*}
\dim H^*_c({\widetilde{\Sigma}''}_v)_{\theta^{-1},\theta}=\dim 
H^*_c(({\widetilde{\Sigma}''}_v)^{H^{\circ}})_{\theta^{-1},\theta}.
\end{equation*}
The Lang--Steinberg theorem implies that both the first and the second projections of $H^{\circ}$ to $T_1$ are surjective. Therefore $(x,x',u',u^{-},\tau,u)\in({\widetilde{\Sigma}''}_v)^{H^{\circ}}$ only if $x=x'=u'=u^{-}=u=1$. Thus $({\widetilde{\Sigma}''}_v)^{H^{\circ}}=\{(1,1,1,1,\tau,1)\mid F(\hat{v}\tau)=\hat{v}\tau\}^{H^{\circ}}$. The set $(\hat{v}T)^F$ is empty unless $\hat{v}^{-1}F(\hat{v})\in T$, in which case $\{(1,1,1,1,\tau,1)\mid F(\hat{v}\tau)=\hat{v}\tau\}$ is actually stable under the action of $H$, so it is also stable under the action of $H^{\circ}$. We only need to treat the non-empty case. As a finite set $(\hat{v}T)^F$ admits only the trivial action of the connected non-trivial group $H^{\circ}$, thus 
\begin{equation*}
({\widetilde{\Sigma}''}_v)^{H^{\circ}}=\{(1,1,1,1,\tau,1)\mid 
F(\hat{v}\tau)=\hat{v}\tau\}^{H^{\circ}}\cong(\hat{v}T)^F.
\end{equation*}
Therefore 
$H^*_c({\widetilde{\Sigma}''}_v)=\overline{\mathbb{Q}}_{\ell}[(\hat{v}T)^F]$, 
on which $T^F\times T^F$ acts via
$(t,t')\colon \hat{v}\tau\mapsto \hat{v}(t^{\hat{v}})^{-1}\tau t'$; note that 
this is the regular representation of both the left $T^F$ and the right $T^F$ 
in $T^F\times T^F$. In particular, the irreducible constituents of 
$H^*_c({\widetilde{\Sigma}''}_v)$ are of the form 
$H^*_c({\widetilde{\Sigma}''}_v)_{(\phi^{\hat{v}})^{-1},\phi}$, where $\phi$ 
runs over $\widehat{T^F}$. Hence  
$H^*_c({\widetilde{\Sigma}''}_v)_{\theta^{-1},\theta}$ is non-zero if and only 
if $\theta^{\hat{v}}=\theta$. As $\theta$ is assumed to be in general position, 
this is equivalent to $v=1$. For $v=1$, we have $\dim 
H^*_c({\widetilde{\Sigma}''}_1)_{\theta^{-1},\theta}=1$ for any 
$\theta\in\widehat{T^F}$, since $|\widehat{T^F}|=|T^F|$. This proves (a).
\end{proof}

To show (b), we use a general homotopy result from \cite{DL1976}:

\begin{lemm}\label{homotopy property}
Let $H$ be a connected algebraic group over $k$, and $Y$ a separated scheme of 
finite type over $k$. Suppose there is a morphism $f\colon H\times Y\rightarrow 
Y$ such that $f(1,-)$ is the identity map and $(h,y)\mapsto (h,f(h,y))$ is an 
automorphism on $H\times Y$. Then for any $h\in H$, the induced endomorphism of 
$f(h,-)$ on $H^i_c(Y,\overline{\mathbb{Q}}_{\ell})$ is the identity map.
\end{lemm}
\begin{proof}
The same argument as in \cite[p.~136]{DL1976} works here.
\end{proof}

To proceed with the proof of the theorem, we need a variant of \cite[Lemma 1.7]{Lusztig2004RepsFinRings}. For general linear groups this can be done in an ad hoc way explicitly (see \cite{ZheChen_PhDthesis}); for general reductive groups we prove the below lemma.

\begin{defi}\label{nottation fixing}
Here we fix several pieces of notation:
\begin{itemize}
\item[(1)] Suppose $\Phi^{-}$ (negative roots of $\mathbf{T}$) is equipped with a total order. For $z\in U^{-}$ and $\beta\in\Phi^{-}$, define $x^{z}_{\beta}\in U_{\beta}=\mathcal{F}\mathbf{U}_{\beta}$ by the decomposition $z=\prod_{\beta\in \Phi^{-}}x^{z}_{\beta}$, where the product is with respect to the following order: If $\mathrm{ht}(\beta)<\mathrm{ht}(\beta')$, then $x^{z}_{\beta}$ is to the left of $x^{z}_{\beta'}$; and if $\mathrm{ht}(\beta)=\mathrm{ht}(\beta')$ and $\beta<\beta'$, then $x^{z}_{\beta}$ is to the left of $x^{z}_{\beta'}$. 

\item[(2)] For a fixed $\alpha\in\Phi^{+}$ and $i\in\{0,...,l-1\}$, denote by $Z^{\alpha}(i)\subseteq U^{-}$ the subvariety consisting of all $z$ such that:
\begin{itemize}
\item[i.] $z\in (U^{-})^i\setminus (U^{-})^{i+1}$;
\item[ii.] $x_{-\alpha}^z\neq1$;
\item[iii.] $x^{z}_{\beta}=1$ for $\forall\beta\in\Phi^{-}$ such that $\mathrm{ht}(\beta)<\mathrm{ht}(-\alpha)$;
\item[iv.] $x^{z}_{\beta}=1$ for $\forall\beta\in\Phi^{-}$ such that $\mathrm{ht}(\beta)=\mathrm{ht}(-\alpha)$ and $\beta<-\alpha$.
\end{itemize}
\end{itemize}
\end{defi}

Recall that $\mathcal{T}^{\alpha}:=(\mathcal{F}\mathbf{T}^{\alpha})^{r-1}$ is a 
$1$-dimensional affine space.

\begin{lemm}\label{main lemma for even level}
Suppose $\alpha\in\Phi^{+}$ and $i\in\{0,...,l-1\}$. Then for $z\in 
Z^{\alpha}(i)$ and $\xi\in U_{\alpha}^{r-i-1}$, one has
\begin{equation*}
[\xi,z]:=\xi z\xi^{-1}z^{-1}=\tau_{\xi,z}\omega_{\xi,z},
\end{equation*}
where $\tau_{\xi,z}\in \mathcal{T}^{\alpha}$ and $\omega_{\xi,z}\in 
{(U^{-})^{r-1}}$ are uniquely determined. Moreover,
\begin{equation*}
U_{\alpha}^{r-i-1}\longrightarrow\mathcal{T}^{\alpha},\qquad 
\xi\longmapsto\tau_{\xi,z}
\end{equation*}
is a surjective morphism admitting a section $\Psi^{\alpha}_z$ such that $\Psi^{\alpha}_z(1)=1$ and such that the map  
\begin{equation*}
Z^{\alpha}(i)\times \mathcal{T}^{\alpha}\longrightarrow U^{r-i-1}_{\alpha},\qquad(z,\tau)\longmapsto \Psi_z^{\alpha}(\tau)
\end{equation*}
is a morphism.
\end{lemm}

\begin{proof}
Write $z=x^z_{-\alpha}z'$, then 
\begin{equation}\label{temp1}
[\xi,z]=\xi x^z_{-\alpha}z'\xi^{-1}z'^{-1}(x^z_{-\alpha})^{-1}=[\xi,x^z_{-\alpha}]\cdot{^{x^z_{-\alpha}}[\xi,z']}.
\end{equation}
We need to determine $[\xi,x^z_{-\alpha}]$ and ${^{x^z_{-\alpha}}[\xi,z']}$ separately.

\vspace{2mm} Following the notation in \cite[XX]{SGA3} we write 
$p_{\beta}\colon 
(\mathbb{G}_a)_{\mathcal{O}^{\mathrm{ur}}_r}\cong\mathbf{U}_{\beta}$ for 
every $\beta\in\Phi$ (and we use the same notation for the isomorphism 
$\mathcal{F}(\mathbb{G}_a)_{\mathcal{O}^{\mathrm{ur}}_r}\cong U_{\beta}$ 
induced by $p_{\beta}$  via the Greenberg functor). Then for some 
$a\in\mathbb{G}_m(\mathcal{O}^{\mathrm{ur}}_r)$ we have
\begin{equation}\label{product of opposite root subgroup}
p_{\alpha}(x)p_{-\alpha}(y)=p_{-\alpha}(\frac{y}{1+axy})\check{\alpha}(1+axy)p_{
\alpha}(\frac{x}{1+axy}),
\end{equation}
for $\forall x,y\in \mathbb{G}_a(\mathcal{O}^{\mathrm{ur}}_r)$; see \cite[XX 2.2]{SGA3}. Let $x,y$ be such that $p_{\alpha}(x)=\xi$ and $p_{-\alpha}(y)=x^z_{-\alpha}$ (note that in our case $x^2=0$ and $(1+axy)^{-1}=1-axy$). By applying (\ref{product of opposite root subgroup}) to the commutator $[p_{\alpha}(x),p_{-\alpha}(y)]=p_{\alpha}(x)p_{-\alpha}(y)p_{\alpha}(-x)p_{-\alpha}(-y)$, we see that
\begin{equation}\label{temp2}
\begin{split}
[\xi,x^z_{-\alpha}]&=p_{\alpha}(x)p_{-\alpha}(y)p_{\alpha}(-x)p_{-\alpha}(-y)\\
&=p_{\alpha}(x)p_{\alpha}(\frac{-x}{1-axy})\check{\alpha}(1+axy)p_{-\alpha}(\frac{y}{1-axy})p_{-\alpha}(-y)\\
&=\check{\alpha}(1+axy)p_{-\alpha}(axy^2).
\end{split}
\end{equation}
Note that since $\xi\in G^{r-i-1}$ and $x_{-\alpha}^z\in G^i$, we have $p_{-\alpha}(axy^2)\in U_{-\alpha}^{r-1}$. In the below one will see that $\check{\alpha}(1+axy)$ is the required $\tau_{\xi,z}$.

\vspace{2mm} Now turn to $[\xi,z']$; we want to show that $[\xi,z']\in (U^{-})^{r-1}$. Let us do this by induction on $\#\{\beta\in\Phi^{-}\mid x^{z'}_{\beta}\neq1\}$. 

\vspace{2mm} When $\#\{\beta\in\Phi^{-}\mid x^{z'}_{\beta}\neq1\}=1$, we have $z'=x^{z'}_{\beta}$ for some $\beta\in\Phi^{-}$, so by the Chevalley commutator formula (\cite[Lemma 2.9 (b)]{Sta2009Unramified}) we have 
\begin{equation*}
[\xi,z']\in\prod_{j,j'\geq1,\  j\beta+j'\alpha\in\Phi}U^{r-1}_{j\beta+j'\alpha}.
\end{equation*}
By basic properties of root systems, if $j\beta+j'\alpha\in\Phi^{+}$ for some 
$j,j'$, then $\beta+\alpha\in\Phi^{+}$. This implies $\mathrm{ht}(-\alpha)>\mathrm{ht}(\beta)$, which is a contradiction to our assumption on $z$, so $[\xi,z']\in (U^-)^{r-1}$ in this case. 

\vspace{2mm} Suppose now that $[\xi,z']\in (U^-)^{r-1}$ whenever $\#\{\beta\in\Phi^{-}\mid x^{z'}_{\beta}\neq1\}\leq N$. Then in case $\#\{\beta\in\Phi^{-}\mid x^{z'}_{\beta}\neq1\}=N+1$, we can decompose the product $z'=\prod_{\beta\in\Phi^{-}}x^{z'}_{\beta}=z'_1z'_2$ such that both $[\xi,z'_1]$ and $[\xi,z'_2]$ are in $(U^-)^{r-1}$. Note that 
\begin{equation*}
[\xi,z']=[\xi,z'_1]\cdot{^{z'_1}[\xi,z'_2]}.
\end{equation*}
Since $z'_1\in U^{-}$, we see that ${^{z'_1}[\xi,z'_2]}\in (U^{-})^{r-1}$, thus $[\xi,z']\in (U^-)^{r-1}$ also for $\#\{\beta\in\Phi^{-}\mid x^{z'}_{\beta}\neq1\}=N+1$. So by induction principle, $[\xi,z']\in (U^-)^{r-1}$ in general.

\vspace{2mm} By \eqref{temp1} and \eqref{temp2} we have
\begin{equation*}
[\xi,z]=[\xi,x^z_{-\alpha}]\cdot{^{x^z_{-\alpha}}[\xi,z']}=\check{\alpha}(1+axy)\cdot p_{-\alpha}(axy^2)\cdot{^{x^z_{-\alpha}}[\xi,z']}.
\end{equation*}
From this expression, put
\begin{equation*}
\tau_{\xi,z} =\check{\alpha}(1+axy)
\end{equation*}
and 
\begin{equation*}
\omega_{\xi,z} =p_{-\alpha}(axy^2)\cdot{^{x^z_{-\alpha}}[\xi,z']}.
\end{equation*}
Note that $\tau_{\xi,z}\in\mathcal{T}^{\alpha}$ and $\omega_{\xi,z}\in (U^-)^{r-1}$ (since $[\xi,z']\in (U^{-})^{r-1}$). The elements $\tau_{\xi,z}$ and $\omega_{\xi,z}$ are uniquely determined because of the Iwahori decomposition.

\vspace{2mm} Now, as $\tau_{\xi,z}$ is defined to be $\check{\alpha}(1+ap^{-1}_{\alpha}(\xi)p^{-1}_{-\alpha}(x^z_{-\alpha}))$, the map $\xi\mapsto\tau_{\xi,z}$, whose target is a connected $1$-dimensional algebraic group, is a surjective algebraic group morphism (note that $z\mapsto x_{-\alpha}^z$ is a projection, hence a morphism). The section morphism $\Psi^{\alpha}_z$ can be defined in the following way: The isomorphism of additive groups 
\begin{equation*}
(\pi^i)\cong \mathcal{O}^{\mathrm{ur}}_{r-i},\quad 
\pi^ia+(\pi^r)\longmapsto a+(\pi^{r-i})
\end{equation*}
induces an isomorphism of affine spaces (by the Greenberg functor)
\begin{equation*}
\mu_i\colon 
(\mathcal{F}(\mathbb{G}_{a})_{\mathcal{O}^{\mathrm{ur}}_r}
)^i\longrightarrow(\mathcal{F}(\mathbb{G}_{a})_{\mathcal{O}^{\mathrm{ur}}_r}
)_{r-i}.
\end{equation*}
Note that this isomorphism depends on the choice of $\pi$. Meanwhile, let
\begin{equation*}
\mu^i\colon (\mathcal{F}(\mathbb{G}_{a})_{\mathcal{O}^{\mathrm{ur}}_r})_{r-i}\cong\mathcal{F}(\mathbb{G}_{a})_{\mathcal{O}^{\mathrm{ur}}_r}/(\mathcal{F}(\mathbb{G}_{a})_{\mathcal{O}^{\mathrm{ur}}_r})^{r-i}\longrightarrow\mathcal{F}(\mathbb{G}_{a})_{\mathcal{O}^{\mathrm{ur}}_r}
\end{equation*}
be a section morphism to the quotient morphism such that $\mu^i(0)=0$ ($\mu^i$ exists because $\mathcal{F}(\mathbb{G}_{a})_{\mathcal{O}^{\mathrm{ur}}_r}$ is an affine space). For $\tau\in\mathcal{T}^{\alpha}$ we put
\begin{equation*}
\Psi_{z}^{\alpha}(\tau):=p_{\alpha}\left(a^{-1}\cdot\mu^{i}\left(  \mu_i\left( \check{\alpha}^{-1}(\tau)-1  \right)\cdot \mu_i\left(  p_{-\alpha}^{-1}(x^z_{-\alpha}) \right)^{-1}   \right) \right).
\end{equation*}
Here $\check{\alpha}^{-1}$ is defined on $\mathcal{T}^{\alpha}=(\mathcal{F}\mathbf{T}^{\alpha})^{r-1}\cong(\mathcal{F}(\mathbb{G}_{m})_{\mathcal{O}^{\mathrm{ur}}_r})^{r-1}$ as the inverse to $\check{\alpha}$, and we view $\check{\alpha}^{-1}(\tau)$ as an element in $\mathcal{F}(\mathbb{G}_{a})_{\mathcal{O}^{\mathrm{ur}}_r}$ by the natural open immersion $(\mathbb{G}_{m})_{\mathcal{O}^{\mathrm{ur}}_r}\rightarrow(\mathbb{G}_{a})_{\mathcal{O}^{\mathrm{ur}}_r}$, so the minus operation $ \check{\alpha}^{-1}(\tau)-1$ is well-defined. On the other hand, by our assumption on $z$ (see Definition \ref{nottation fixing} (2) i), $\mu_i\left( p_{-\alpha}^{-1}(x^z_{-\alpha}) \right)$ is an element in $\mathcal{F}(\mathbb{G}_{m})_{\mathcal{O}^{\mathrm{ur}}_{r-i}}$, so its multiplicative inverse exists. Moreover, the product operation ``$\cdot$'' is by viewing $(\mathbb{G}_{a})_{\mathcal{O}^{\mathrm{ur}}_r}$ (resp.\ $\mathcal{F}(\mathbb{G}_{a})_{\mathcal{O}^{\mathrm{ur}}_r}$) as a ring scheme (resp.\ $k$-ring variety). Thus $\Psi_z^{\alpha}$ is well-defined as a morphism. 

\vspace{2mm} Finally, by the definition of $\mu_i$ and $\mu^i$, for 
$\tau\in\mathcal{T}^{\alpha}(k)$ we have
\begin{equation*}
\begin{split}
\tau_{\Psi_z^{\alpha}(\tau),z}&=\check{\alpha}\left(1+ap_{-\alpha}^{-1}(\Psi_z^{
\alpha}(\tau))p^{-1}_{-\alpha}(x_{-\alpha}^z)\right)\\
&=\check{\alpha}\left(1+\mu^{i}\left(  \mu_i\left(  \check{\alpha}^{-1}(\tau)-1 
 \right)\cdot \mu_i\left(  p_{-\alpha}^{-1}(x^z_{-\alpha}) \right)^{-1}  
\right) \cdot p^{-1}_{-\alpha}(x_{-\alpha}^z) \right)\\
&=\check{\alpha}\left(1+\pi^i\cdot\mu^i\mu_i(\check{\alpha}^{-1}
(\tau)-1)\right)=\tau
\end{split}
\end{equation*}
(for the last equality, note that $\check{\alpha}^{-1}(\tau)$ is of the form $1+s\pi^{r-1}$ for some $s\in\mathcal{O}^{\mathrm{ur}}_r$, as an element in $\mathbb{G}_m(\mathcal{O}^{\mathrm{ur}}_r)$), thus $\tau\mapsto\Psi_z^{\alpha}(\tau)\mapsto\tau_{\Psi_z^{\alpha}(\tau),z}$ is the identity map on the $k$-points $\mathcal{T}^{\alpha}(k)$ of the $1$-dimensional affine space $\mathcal{T}^{\alpha}\cong\mathbb{A}^1_k$, hence it is the identity morphism. So $\Psi_{z}^{\alpha}$ is a section to $\xi\mapsto\tau_{\xi,z}$, and the other assertions in the lemma follow from its definition.
\end{proof}

\begin{lemm}\label{lemm b}
(b) is true.
\end{lemm}
\begin{proof}
By the changes of variables $\hat{v}\tau\hat{v}^{-1}\mapsto\tau$, 
$\tau^{-1}u^{-}\tau\mapsto u^{-}$, and $\tau^{-1}u'\tau\mapsto u'$, we can 
rewrite $\Sigma'_v$ as
\begin{equation*}
\widetilde{\Sigma}'_v:=\{(x,x',u',u^{-},\tau,u)\in U^{\pm}\times FU\times 
\mathcal{Z}_v\setminus\mathcal{Z}_v(l)\times T\times U\mid xF(\tau 
u'u^{-}\hat{v})=\tau u'u^{-}\hat{v}ux'\},
\end{equation*}
on which $(t,t')\in T^F\times T^F$ acts by sending $(x,x',u',u^{-},\tau,u)$ to
\begin{equation*} 
(t^{-1}xt,{t'}^{-1}x't',({t'}^{\hat{v}})^{-1}u'(t')^{\hat{v}},({t'}^{\hat{v}})^{
-1}u^{-}(t')^{\hat{v}},t^{-1}\tau (t')^{\hat{v}},{t'}^{-1}ut').
\end{equation*}
To show $\dim H^*_c(\widetilde{\Sigma}'_v)_{\theta^{-1},\theta}=0$, it suffices 
to show 
\begin{equation*}
\dim H^*_c(\widetilde{\Sigma}'_v)_{\theta^{-1}|_{(T^{r-1})^F}}=0,
\end{equation*}
for the subgroup $(T^{r-1})^F=(T^{r-1})^F\times 1\subseteq T^F\times T^F$. Note 
that the $(T^{r-1})^F$-action on $\widetilde{\Sigma}'_v$ is given by 
\begin{equation*}
t\colon (x,x',u',u^{-},\tau,u)\mapsto (x,x',u',u^{-},t^{-1}\tau,u).
\end{equation*}

\vspace{2mm} Recall that we fixed an order on $\Phi^{-}$. For 
$\beta\in\Phi^{-}$, let $F(\beta)\in\Phi$ be the root defined by 
$F({U})_{F(\beta)}=F({U}_{\beta})$, then the order on $\Phi^-$ produces an 
order on $F(\Phi^-)$; similarly we can define $F$ on $\Phi^{+}$, and hence get 
a bijection on $\Phi=\Phi^-\sqcup\Phi^+=F(\Phi^-)\sqcup F(\Phi^+)$, and then a 
bijection on $\{{U}_{\beta}\}_{\beta\in \Phi}$; it is clear that
$F(-\alpha)=-F(\alpha)$ for any $\alpha\in\Phi$. Let $\mathcal{Z}_v^{\beta}(i)$ 
be the subvariety of $\mathcal{Z}_v(i)\setminus\mathcal{Z}_v(i+1)$ consisted of 
$(u',u^{-})$ such that, in the decomposition 
$F(z):=F(\hat{v}^{-1}u'u^{-}\hat{v})=\prod_{\beta'\in{F(\Phi^{-})}}x^{F(z)}_{
\beta'}$ one has: $x_{\beta'}^{F(z)}=1$ whenever 
$\mathrm{ht}(\beta')<\mathrm{ht}(F(\beta))$, $x_{\beta'}^{F(z)}=1$ whenever 
$\mathrm{ht}(\beta')=\mathrm{ht}(F(\beta))$ and $\beta'<F(\beta)$, and 
$x_{F(\beta)}^{F(z)}\neq1$ (compare the conditions in Definition \ref{nottation 
fixing} (2) by formally replacing $\alpha$ by 
$-F(\beta)$ and $\Phi^{-}$ by $F(\Phi^{-})$). We then obtain a finite partition
\begin{equation*}
\mathcal{Z}_v\setminus\mathcal{Z}_v(l)=\coprod_{i=0}^{l-1}\coprod_{\beta\in\Phi^
{-}}\mathcal{Z}_v^{\beta}(i).
\end{equation*} 
And hence a partition of $\widetilde{\Sigma}'_v$ into locally closed 
subvarieties
\begin{equation*}
\widetilde{\Sigma}'_v=\coprod_{i=0}^{l-1}\coprod_{\beta\in\Phi^{-}}
\Sigma_v^{\beta}(i),
\end{equation*} 
where
\begin{equation*}
\Sigma_v^{\beta}(i):=\{(x,x',u',u^{-},\tau,u)\in U^{\pm}\times FU\times 
\mathcal{Z}_v^{\beta}(i)\times T\times U\mid xF(\tau u'u^{-}\hat{v})=\tau 
u'u^{-}\hat{v}ux'\}.
\end{equation*} 
Each subvariety $\Sigma_v^{\beta}(i)$ inherits the $(T^{r-1})^F$-action:
\begin{equation*}
t\colon (x,x',u',u^{-},\tau,u)\mapsto (x,x',u',u^{-},t^{-1}\tau,u),
\end{equation*}
so it suffices to show: 
\begin{equation*}
H^*_c(\Sigma_v^{\beta}(i))_{\theta^{-1}|_{(T^{r-1})^F}}=0
\end{equation*} 
for every $i\in\{0,...,l-1\}$ and every $\beta\in \Phi^{-}$.

\vspace{2mm} From now on we fix an $\alpha\in\Phi^{+}$. Consider the closed subgroup 
\begin{equation*}
H:=\{t\in T^{r-1}\mid F(\hat{v})^{-1}F(t)t^{-1}F(\hat{v})\in \mathcal{T}^{F(\alpha)}\}
\end{equation*}
of $T^{r-1}$. For any $t\in H$, define $g_t\colon FU\rightarrow FU$ by
\begin{equation*}
g_t\colon x'\mapsto x'\cdot\Psi_{F(z)}^{F(\alpha)}\left(F(\hat{v})^{-1}F(t^{-1})tF(\hat{v})\right)^{-1}
\end{equation*}
with the parameter $z:=\hat{v}^{-1}u'u^{-}\hat{v}$, where 
$(u',u^{-})\in\mathcal{Z}_v^{-\alpha}(i)$. This is well-defined because $F(z)$ 
satisfies the conditions in Lemma \ref{main lemma for even level}, with respect 
to $F(U^{-})$ and $F(\Phi^{-})$. Note that if $F(t)=t$, then $g_t(x')=x'$.

\vspace{2mm} Moreover, for any $t\in H$, define the morphism $f_t\colon 
U^{\pm}\rightarrow U^{\pm}$ by
\begin{equation*}
f_t\colon x\mapsto x\cdot {^{F(\tau)}\left( t^{-1} \cdot {^{F(\hat{v}z)}\left( x'^{-1}g_t(x')\right)}  F(t) \right)},
\end{equation*}
with the parameters $x'\in FU$, $\tau\in T$, and $z=\hat{v}^{-1}u'u^{-}\hat{v}$ 
(where $(u',u^{-})\in\mathcal{Z}_v^{-\alpha}(i)$, as for $g_t$). To see this is 
well-defined one needs to check the right hand side is in $U^{\pm}$: By the 
definition of $\Psi_{F(z)}^{F(\alpha)}$ and the first assertion of Lemma 
\ref{main lemma for even level} we see
\begin{equation*}
F(z) x'^{-1}g_t(x') F(z^{-1})=\Psi_{F(z)}^{F(\alpha)}\left(F(\hat{v})^{-1}F(t^{-1})tF(\hat{v})\right)^{-1}\cdot F(\hat{v})^{-1}F(t^{-1})tF(\hat{v})\cdot\omega
\end{equation*}
for some $\omega\in{U^{r-1}}$. Hence by definition of $f_t$ we get
\begin{equation*}
(x^{-1}f_t(x))^{F(\tau)}=({^{F(\hat{v})}\Psi})^t   \cdot ({^{F(\hat{v})}\omega})^{F(t)} \in\prod_{\beta\in\Phi}{U^{r-i-1}_{\beta}}\subseteq U^{\pm},
\end{equation*}
where $\Psi:=\Psi_{F(z)}^{F(\alpha)}(F(\hat{v})^{-1}F(t^{-1})tF(\hat{v}))^{-1}$. Thus $x^{-1}f_t(x)\in U^{\pm}$, and $f_t$ is therefore well-defined. Moreover, if $F(t)=t$, then $f_t(x)=x$.

\vspace{2mm} For any $t\in H$, the above preparations on $f_t$ and $g_t$ allow us to define the following automorphism of $\Sigma^{-\alpha}_v(i)$:
\begin{equation*}
h_t\colon (x,x',u',u^{-},\tau,u)\mapsto (f_t(x),g_t(x'),u',u^{-},t^{-1}\tau,u),
\end{equation*}
where the involved parameter $z$ is $\hat{v}^{-1}u'u^{-}\hat{v}$. To see this is well-defined, one needs to show the right hand side satisfies the defining equation of $\Sigma^{-\alpha}_v(i)$, in other words, satisfies 
\begin{equation*}
f_t(x)F(t^{-1}\tau u'u^{-}\hat{v})=t^{-1}\tau u'u^{-}\hat{v}ug_t(x');
\end{equation*}
this can be seen by just expanding the definition of $f_t$: (note that $t\in T^{r-1}$ commutes with $x\in U^{\pm}$, and $xF(\tau u'u^{-}\hat{v})=\tau u'u^{-}\hat{v}ux'$)
\begin{equation*}
\begin{split}
f_t(x)F(t^{-1}\tau u'u^{-}\hat{v})&=x\cdot {^{F(\tau)}\left( t^{-1} \cdot {^{F(\hat{v}z)}\left( x'^{-1}g_t(x')\right)} F(t) \right)}\cdot F(t^{-1}\tau u'u^{-}\hat{v})\\
&=t^{-1}xF(\tau u'u^{-}\hat{v})x'^{-1}g_t(x')\\
&=t^{-1}\tau u'u^{-}\hat{v}ug_t(x').
\end{split}
\end{equation*}
Moreover, it is clear that in the case $F(t)=t$, the automorphism $h_t$ coincides with the $(T^{r-1})^F$-action, so by Lemma \ref{homotopy property}, the induced endomorphism of $h_t$ on $H^*_c(\Sigma^{-\alpha}_v(i))$ is the identity map for any $t$ in the identity component $H^{\circ}$ of $H$.

\vspace{2mm}Let $a\geq1$ be an integer such that $F^a(F(\hat{v})\mathcal{T}^{F(\alpha)}F(\hat{v})^{-1})=F(\hat{v})\mathcal{T}^{F(\alpha)}F(\hat{v})^{-1}$, then the image of the norm map $N^{F^a}_F(t)=t\cdot F(t)\cdots F^{a-1}(t)$ on $F(\hat{v})\mathcal{T}^{F(\alpha)}F(\hat{v})^{-1}$ is a connected subgroup of $H$, hence contained in $H^{\circ}$. Moreover $N^{F^a}_F((F(\hat{v})\mathcal{T}^{F(\alpha)}F(\hat{v})^{-1})^{F^a})\subseteq (T^{r-1})^F\cap H^{\circ}$. Thus, as $\theta$ is regular, 
\begin{equation*}
H^*_c(\Sigma^{\beta}_v(i))_{\theta^{-1}\big|_{ N^{F^a}_F\left(\left(F(\hat{v})\mathcal{T}^{F(\alpha)}F(\hat{v})^{-1}\right)^{F^a}\right)}}=0.
\end{equation*}
Therefore $H^*_c(\Sigma^{\beta}_v(i))_{\theta^{-1}|_{(T^{r-1})^F}}=0$, which proves $(b)$.
\end{proof}

By Lemma \ref{lemm a} and Lemma \ref{lemm b}, $\dim H^*_c(\Sigma)_{\theta^{-1},\theta}=\dim H^*_c(\Sigma''_1)_{\theta^{-1},\theta}=1$. Thus we complete the proof of the theorem.
\end{proof}

\begin{coro}\label{Main result at even levels}
Let $r=2l$, and suppose $\theta\in\widehat{T^F}$ is a generic character; denote by $\widetilde{\theta}$ the trivial lift of $\theta$ to $(TU^{\pm})^F=(TG^l)^F$. Then we have $R_{T}^{\theta}\cong\mathrm{Ind}_{(TU^{\pm})^F}^{G^F}\widetilde{\theta}$, and they are irreducible representations of dimension $|G_l^F|/|T_l^F|$.
\end{coro}
\begin{proof}
As $\theta$ is generic, $R_{T}^{\theta}$ is irreducible by Theorem \ref{Irr DL}, and $\mathrm{Ind}_{(TU^{\pm})^F}^{G^F}\widetilde{\theta}$ is irreducible by Clifford theory. So the result follows from Theorem \ref{multiplicity one, even levels}.
\end{proof}

\section{An application to finite Lie algebras}\label{finite lie alg}

In this last section we assume $\mathcal{O}=\mathbb{F}_q[[\pi]]$ and $r=2$. Note 
that the kernel group $G^{1}$ is isomorphic to the additive group of the Lie 
algebra $\mathfrak{g}$ of $G_1$, and the adjoint action of $G_1^F$ on 
$\mathfrak{g}^F$ is the conjugation action under this isomorphism. Since 
$T^F\cong T_1^F\times (T^1)^F$, any character $\theta^1$ of 
${\mathfrak{t}^F}\cong(T^{1})^F$ extends (trivially) to a character ${\theta}$ 
of $T^F$. Thus, by viewing $R_{T,U}^{{\theta}}$ as a 
$\mathfrak{g}^F\cong(G^{1})^F$-module, we can view 
$R_{\mathfrak{t},\mathfrak{u}}^{\theta^1}:=R_{T,U}^{{\theta}}$ as a 
Deligne--Lusztig theory for the finite Lie algebra $\mathfrak{g}^F$ (here 
$\mathfrak{u}$ is the Lie algebra of $U_1$). 

\vspace{2mm} An \emph{invariant character} of $\mathfrak{g}^F$ is a $\overline{\mathbb{Q}}_{\ell}$-character of the finite abelian group $\mathfrak{g}^F$ that is invariant under the adjoint action of $G_1^F$, and it is said to be irreducible if it is not the sum of two invariant characters (these functions have interesting relations with character sheaves; see e.g.\ \cite{Lusztig_1987_fourier_Lie_alg} and \cite{Let2005book}). Letellier studied this construction in \cite{Let09CharRedLieAlg}, where he compared this construction with a different construction he considered earlier in \cite{Let2005book}, and made a conjecture that every irreducible invariant character $\Psi$ of $\mathfrak{g}^F$ ``appear'' in some $R_{\mathfrak{t},\mathfrak{u}}^{\theta^1}$ in the sense that 
\begin{equation*}
(\Psi,R_{\mathfrak{t},\mathfrak{u}}^{\theta^1})_{\mathfrak{g}^F}:=\frac{1}{|G_1^F|}\sum_{g\in\mathfrak{g} ^F}\Psi(g)R_{\mathfrak{t},\mathfrak{u}}^{\theta^1}(-g)\neq0
\end{equation*} 
(note that the bracket $(,)$ is different from the usual inner product 
$\langle,\rangle$). Letellier's result shows this conjecture is true for 
$\mathrm{GL}_2$ with the assumption that $\mathrm{char}(\mathbb{F}_q)>3$. Here 
as a simple application of our main result, we remove this assumption.

\begin{prop}
Along with the above notation, if $\mathbb{G}=\mathrm{GL}_2$ or $\mathrm{GL}_3$, then for any irreducible invariant character $\Psi$ of $\mathfrak{g}^F$, we have
\begin{equation*}
(\Psi,R^{\theta^1}_{\mathfrak{t},\mathfrak{u}})_{\mathfrak{g}^F}\neq0,
\end{equation*}
for some $R^{\theta^1}_{\mathfrak{t},\mathfrak{u}}$.
\end{prop}
\begin{proof}
Firstly note that $(\Psi,R^{\theta^1}_{\mathfrak{t},\mathfrak{u}})_{\mathfrak{g}^F}\neq0$ if and only if $\langle\Psi,R^{\theta^1}_{\mathfrak{t},\mathfrak{u}}\rangle_{(G^1)^F}\neq0$ (these are two different brackets). Also note that a $\mathfrak{g}^F$-representation is invariant if and only if it is $G^F$-invariant as a $(G^{1})^F$-representation, so we can focus on characters of the group $(G^{1})^F$. Suppose $\chi$ is an irreducible character of $(G^{1})^F$, then 
\begin{equation*}
\chi^{O}:=\bigoplus_{s\in G^F/\mathrm{Stab}_{G^F}(\chi)}\chi^s
\end{equation*}
is an invariant character of $(G^{1})^F$, and any invariant character containing $\chi$ contains $\chi^{O}$ (so $\chi^{O}$ is the unique irreducible invariant character containing $\chi$). On the other hand, any $G^F$-module is an invariant $(G^{1})^F$-module, thus we only need to show any irreducible character $\chi$ of $(G^{1})^F$ is ``contained'' in some $R_{\mathfrak{t},\mathfrak{u}}^{\theta^1}$ in the sense that $\langle\chi,R^{\theta^1}_{\mathfrak{t},\mathfrak{u}}\rangle_{(G^1)^F}\neq0$. 

\vspace{2mm} For $\mathbb{G}=\mathrm{GL}_2$ (resp.\ $\mathrm{GL}_3$), the irreducible characters of $\mathfrak{g}^F$ are of the form $\chi=\psi_{\beta}(-)=\psi(\mathrm{Tr}(\beta\cdot(-)))$, where $\psi$ is some fixed non-trivial $\overline{\mathbb{Q}}_{\ell}$-character of $\mathbb{F}_q$ and $\beta\in M_2(\mathbb{F}_q)$ (resp.\ $\beta\in M_3(\mathbb{F}_q)$). The conjugacy classes of $\beta\in M_2(\mathbb{F}_q)$ are of the following two types:
\begin{itemize}
\item[\rm{(1)}]
$\begin{bmatrix}
a & * \\
0 & b \\
\end{bmatrix}$, where $*$ is $0$ or $1$;

\item[\rm{(2)}]
$\begin{bmatrix}
0 & 1 \\
-\Delta & s \\
\end{bmatrix}$, where $x^2-sx+\Delta$ is irreducible over $\mathbb{F}_q$.
\end{itemize}
And the conjugacy classes of $\beta\in M_3(\mathbb{F}_q)$ are of the following three types:
\begin{itemize}
\item[\rm{(1')}]
$\begin{bmatrix}
a & *_1 & 0\\
0 & b & *_2\\
0 & 0 & c\\
\end{bmatrix}$, where $*_1$ and $*_2$ are $0$ or $1$;

\item[\rm{(2')}]
$\begin{bmatrix}
0 & 1 & 0\\
-\Delta & s & 0\\
0 & 0 & a
\end{bmatrix}$, where $x^2-sx+\Delta$ is irreducible over $\mathbb{F}_q$;

\item[\rm{(2'')}] $N$, where $\det(x\cdot I-N)$ is irreducible over $\mathbb{F}_q$.
\end{itemize}

\vspace{2mm} For types (1) and (1'), the corresponding $\chi=\psi_{\beta}$ is trivial on the rational points of the Lie algebra of the unipotent radical $\mathbf{U}_0$ of some rational Borel subgroup $\mathbf{B}_0$. Let $\mathbf{T}=\mathbf{T}_0$ be a rational maximal torus contained in $\mathbf{B}_0$, and following the previous notation we denote by $\theta^1$ the restriction of $\chi$ to $\mathfrak{t}^F=(T^1)^F$. Then we have
\begin{equation*}
\left\langle\mathrm{Res}_{\mathfrak{g}^F}^{G^F}\mathrm{Ind}_{B_0^F}^{G^F}\widetilde{{\theta}},\chi\right\rangle_{(G^1)^F}=\sum_{s\in B_0^F\backslash G^F/\mathfrak{g}^F}\left\langle\mathrm{Ind}_{(^s(B_0)^{1})^F}^{\mathfrak{g}^F}\left(\widetilde{{\theta}}^{s^{-1}}|_{(^s(B_0)^{1})^F}\right),\chi\right\rangle_{(G^1)^F}
\end{equation*}
by the Mackey intertwining formula. Note that
\begin{equation*}
\left\langle\mathrm{Ind}_{(^sB_0^{1})^F}^{\mathfrak{g}^F}\left(\widetilde{{\theta}}^{s^{-1}}|_{(^sB_0^{1})^F}\right),\chi\right\rangle_{(G^1)^F}=\left\langle\left(\widetilde{{\theta}}^{s^{-1}}|_{(^sB_0^{1})^F}\right),\chi|_{(^s(B_0)^{1})^F}\right\rangle_{(^s(B_0)^{1})^F}
\end{equation*}
by the Frobenius reciprocity, which is non-zero in the case $s=1$. Therefore $\chi$ appears in $\mathrm{Ind}_{B_0^F}^{G^F}\widetilde{\theta}=R_{\mathfrak{t},\mathfrak{u}}^{\theta^1}$.

\vspace{2mm} For type (2) (resp. types (2'), and (2'')), the $\beta$ is a semisimple regular element in $M_2(\mathbb{F}_q)$ (resp.\ $M_3(\mathbb{F}_q)$), in particular the corresponding $\theta$ is in general position and $\mathrm{Stab}_{G^F}(\theta|_{(G^l)^F})=(TU^{\pm})^F$. For $\mathrm{GL}_2$ (resp.\ $\mathrm{GL}_3$) conjugate $\beta$ to be a diagonal matrix in $M_2(k)$ (resp.\ $M_3(k)$), and view $T^1$ as the set of diagonal matrices in $M_2(k)$ (resp.\ $M_3(k)$) with Frobenius endomorphism being the canonical one conjugated by an element in the Weyl group, then the same argument of Proposition \ref{coxeter generic} shows $\theta$ is regular. So thanks to Corollary \ref{Main result at even levels} we only need to show $\chi=\psi_{\beta}$ appears in $\mathrm{Ind}_{(TU^{\pm})^F}^{G^F}\widetilde{\theta}$. Actually, again by the Mackey intertwining formula we have
\begin{equation*}
\left\langle\mathrm{Res}_{\mathfrak{g}^F}^{G^F}\mathrm{Ind}_{(TU^{\pm})^F}^{G^F}\widetilde{{\theta}},\chi\right\rangle_{(G^1)^F}=\sum_{s\in (TU^{\pm})^F\backslash G^F/\mathfrak{g}^F}\left\langle\widetilde{{\theta}}^{s^{-1}}|_{\mathfrak{g}^F},\chi\right\rangle_{(G^1)^F},
\end{equation*}
which is non-zero (take $s=1$).
\end{proof}

\bibliographystyle{alpha}
\bibliography{zchenrefs,alex}

\end{document}